\documentclass[10pt,twoside]{amsart}
 \usepackage{amsm ath,amssymb,amsthm,amscd,latexsym,graphicx}
 \usepackage{epigraph}
 \usepackage[OT2,OT1]{fontenc}

\usepackage{tikz}
\usepackage{graphicx}
\usepackage[all]{xy}
\usepackage{blindtext}
\usepackage{titletoc}
\usepackage[shortlabels]{enumitem}
 \font\caps=cmcsc10     % Theorem, Lemma etc
 \font\Caps=cmcsc10 scaled \magstep1 % Title

  \entrymodifiers={!!<0pt,0.7ex>+}

 \makeatletter
 \@specialpagefalse\headheight=8.5pt
 \makeatother

 \def\TSkip{\medskip}
 \newbox\TheTitle{\obeylines\gdef\GetTitle #1
 \ShortTitle #2
 \SubTitle #3
 \Author  #4
 \ShortAuthor #5
 \EndTitle
 {\setbox\TheTitle=\vbox{\baselineskip=20pt\let\par=\cr\obeylines%
 \halign{\centerline{\Caps##}\cr\noalign{\medskip}\cr#1\cr}}%
   \copy\TheTitle\TSkip\TSkip%
 \def\next{#2}\ifx\next\empty\gdef\STitle{#1}\else\gdef\STitle{#2}\fi%
 \def\next{#3}\ifx\next\empty%
 \else\setbox\TheTitle=\vbox{\baselineskip=20pt\let\par=\cr\obeylines%
  \halign{\centerline{\caps##} #3\cr}}\copy\TheTitle\TSkip\TSkip\fi%
 \centerline{\caps #4}\TSkip\TSkip%
 \def\next{#5}\ifx\next\empty\gdef\SAuthor{#4}\else\gdef\SAuthor{#5}\fi%
 \catcode'015=5}}

 \long\def\MSC#1\EndMSC{\def\arg{#1}\ifx\arg\empty\relax\else
  {\par\narrower\noindent%
  2000 Mathematics Subject Classification: #1\par}\fi}

 \long\def\KEY#1\EndKEY{\def\arg{#1}\ifx\arg\empty\relax\else
   {\par\narrower\noindent Keywords and Phrases: #1\par}\fi\TSkip}

 \long\def\DATE#1\EndDATE{\def\arg{#1}\ifx\arg\empty\relax\else
   {\par\narrower\noindent \center{\textit{#1}}\par}\fi\TSkip\TSkip\TSkip}

 \newcommand{\Sym}{\hbox{\ptitgot S}}

 \newcommand{\Ind}{\mathrm{Ind}}

\newcommand{\vers}{\mathrm{vers} }

 \newcommand{\Cor}{\mathrm{Cor}}
  \newcommand{\obs}{\mathrm{Obs}}
    \newcommand{\nat}{\mathrm{nat}}
        \newcommand{\aug}{\mathrm{aug}}

 \renewcommand{\to}{\longrightarrow}

 \newcommand{\A}{\mathbb A}

 \newcommand{\F}{\mathbb{F}}
 \newcommand{\G}{\mathbb{G}}

 \newcommand{\N}{\mathbb{N}}
 \renewcommand{\P}{\mathbb P}
 \newcommand{\Q}{\mathbb{Q}}
 
 \newcommand{\Z}{\mathbb{Z}}

 \newcommand{\Res}{\mathrm{Res}}

 \newcommand{\GL}{\mathbf{GL}}

 \newcommand{\RR}{\mathbf{R}}

 \newcommand{\Aut}{\mathrm{Aut}}

 \newcommand{\End}{\mathrm{End}}
 \newcommand{\Ext}{\mathrm{Ext}}
\newcommand{\Gal}{\mathrm{Gal}}
 
 \newcommand{\Hom}{\mathrm{Hom}}

 \newcommand{\Frob}{\mathrm{Frob}}

\newcommand{\Id}{\mathrm{Id}}

\newcommand{\Ker}{\mathrm{Ker}}
 \newcommand{\pr}{\mathrm{pr}}
 \newcommand{\Spec}{\mathrm{Spec}}

\newcommand{\W}{\mathbf{W}}

\newcommand{\Pic}{\mathrm{Pic}}

\renewcommand{\RR}{\mathbf{R}}
\renewcommand{\Sym}{\mathrm{Sym}}

 \theoremstyle{plain}
 \newtheorem{thm}{Theorem}[section]
   \newtheorem*{thm*}{Theorem}
 \newtheorem{defi}[thm]{Definition}

 \newtheorem{prop}[thm]{Proposition}

 \newtheorem{lem}[thm]{Lemma}
 \newtheorem{coro}[thm]{Corollary}
\newtheorem{conj}[thm]{Conjecture.}
 \newtheorem*{theorem-non}{Bloch-Kato conjecture, an equivalent formulation}
\newtheorem*{thmA*}{Theorem A}
\newtheorem*{thmB*}{Theorem B}
\newtheorem*{thmC*}{Theorem C}
\newtheorem*{thmD*}{Theorem D}

 \theoremstyle{remark}
 \newtheorem{rem}[thm]{Remark}
 \newtheorem{qu}[thm]{Question}

 \newtheorem{exo}[thm]{Exercise}

 \newenvironment{dem}{{\bf Proof.}}{\hfill$\square$}

\newcommand{\fonctionnoname}[4]{\begin{array}{ccc}
 #1 & \longrightarrow & #2 \\
 #3 & \longmapsto & #4 \end{array}}

\begin{document}

\title{Smooth profinite groups, III: the Smoothness Theorem.}

\author{Charles De Clercq and Mathieu Florence}
\address{Charles De Clercq, Équipe Topologie Alg\'ebrique, Laboratoire Analyse, G\'eom\'etrie et Applications, Sorbonne Paris Nord, 93430 Villetaneuse, France.}
\address{Mathieu Florence, Équipe topologie et géométrie algébriques, Institut de Mathématiques de Jussieu, Sorbonne Université, F-75005 Paris, France.}


\begin{abstract}
Let $p$ be a prime. In this article, we prove the Smoothness Theorem, which  asserts that a $(1,1)$-cyclotomic pair is $(n,1)$-cyclotomic, for all $n \geq 1$. In the particular case of Galois cohomology, the Smoothness Theorem provides a new proof of the Norm Residue Isomorphism Theorem, entirely disjoint from motivic cohomology. A byproduct of this approach, is that the latter Theorem follows from mod $p^2$ Kummer theory for fields alone.  We moreover extend it, from absolute Galois groups of fields, to algebraic fundamental groups of (not necessarily smooth, nor proper) curves over algebraically closed fields.
\end{abstract}

\maketitle

\tableofcontents
\newpage
\section{Introduction.}

Let $G$ be a profinite group and $p$ be a prime. A decade ago, we began working on a mathematical project rooted in  the following belief.\\

(B): The  Norm Residue Isomorphism Theorem follows from mod $p^2$ Kummer theory for fields.\\

Recall that the Norm Residue Isomorphism Theorem, formerly the Milnor-Bloch-Kato Conjecture, was proved by Rost, Suslin, Voevodsky and Weibel, by applying motivic cohomology to norm varieties \cite{HW}.

With the two previous articles of this series in hand, we can now show that belief (B) is correct, even beyond Galois cohomology: it is accurate in the broader context of smooth profinite groups, introduced in \cite{DCF0}. 

Using the algebro-geometric tools developed in  \cite{F2} (notably the Uplifting Pattern)  we prove the  \textit{Smoothness Theorem} (Theorem \ref{SmoothTh}), a lifting result for mod $p$ cohomology of smooth profinite groups. It states that a $(1,1)$-cyclotomic pair is $(n,1)$-cyclotomic, for any $n\geq 1$. \\

Denote by $G=\Gal(F_s/F)$ `the' absolute Galois group of a field $F$ of characteristic not $p$. Let  $\Z/p^2(1):=\mu_{p^2},$ the Galois module of $p^2$-th roots of unity. Recall that the pair $(G,\Z/p^2(1))$ is  $(1,1)$-cyclotomic (see \cite{DCF1}). The Smoothness Theorem then applies, yielding the Norm Residue Isomorphism Theorem at very low cost, using its equivalent formulation given by Merkurjev for $p=2$ \cite{Me} and for arbitrary $p$ by Gille \cite{Gi}.

Note that we provided in \cite[§3]{DCF0} several classes of profinite groups of geometric origin, that  fit into $(1,1)$-cyclotomic pairs. These include \'etale fundamental groups of semilocal $\mathbb{Z}[\frac{1}{p}]$-schemes and of   connected curves over  algebraically closed fields. [These curves are assumed to be smooth in \textit{loc. cit.} As explained in the Appendix, this is superfluous.] The Smoothness Theorem therefore extends  the validity of the Norm Residue Isomorphism Theorem, to these  \'etale fundamental groups.

We assume familiarity with the  notions developed in \cite{DCF1} and \cite{F2}, of which notation and conventions are kept. In particular, actions of profinite groups on algebro-geometric structures are naive: they factor through open normal subgroups.
\section{The Smoothness Theorem.}\label{SmoothSection}

Let $G$ be a profinite group. One can search for a condition on $G$, for the natural (reduction) arrow on cohomology algebras \[H^*(G,\Z/p^2) \to H^*(G,\F_p)\] to be surjective. The Theorem below provides an  answer- allowing the necessary presence of a twist on  coefficients. This proves  the main part of the Smoothness Conjecture 14.25 of \cite{DCFOr}, in depth $e=1$.

\begin{thm}(The Smoothness Theorem.)\label{SmoothTh}\\
 Let $G$ be a $(1,1)$-smooth profinite group. Then, for every $n \geq 1$,  $G$ is  $(n,1)$-smooth.\\
Moreover, assume that $(G,\Z/p^2(1))$ is a $(1,1)$-cyclotomic pair. \\Then, for every $n \geq 1$,  $(G,\Z/p^2(1))$ is  $(n,1)$-cyclotomic.
\end{thm}

\begin{rem}
    For finitely-generated \textit{analytic} pro-$p$-groups, the Smoothness Theorem was already known- see \cite{Q}. The nice proof,  goes by elementary and purely group-theoretic considerations. See also \cite{QW}.
\end{rem}

\begin{coro}[{The Norm Residue Isomorphism Theorem}]\quad\\
Let $p$ be a prime and $F$ be a field of characteristic not $p$ with separable closure $F_s$.
\begin{enumerate}
    \item For all $n \geq 1$, the natural homomorphism \[ H^n(\Gal(F_s /F),  \mu_{p^2}^{\otimes n}) \to  H^n(\Gal(F_s /F),  \mu_{p}^{\otimes n})\] is onto. Equivalently, the connecting homomorphism \[H^n(\Gal(F_s /F),  \mu_{p}^{\otimes n}) \to H^{n+1}(\Gal(F_s /F),  \mu_{p}^{\otimes n}) , \] arising from the twisted Kummer sequence \[0 \to  \mu_{p}^{\otimes n} \to \mu_{p^2}^{\otimes n} \to  \mu_{p}^{\otimes n} \to 0 ,\] vanishes.\\
    \item The Norm Residue Isomorphism Theorem holds: the Galois symbol \[ K_n^M(F)/p \to H^n(\Gal(F_s /F),  \mu_{p}^{\otimes n})\] is an isomorphism.
\end{enumerate}
\end{coro}

\begin{dem}
Remember that $(\Gal(F_s /F),\mu_{p^2})$ is a $(1,1)$-cyclotomic pair by  Kummer theory (see \cite{DCFOr}). The first statement is then given by the Smoothness Theorem.\\
We give references for proving $1) \Rightarrow 2)$. For $p=2$ (Milnor's Conjecture), this is precisely the work of \cite{Me}, while for $p$ is arbitrary, one can cast \cite[Theorem 0.2]{Gi}. Both papers are short and self-contained.
\end{dem}
\quad\\

 The rest of this paper is devoted to the proof of the Smoothness Theorem. We first deal with the cyclotomic case,   assuming that the cyclotomic module is trivial. [For $G=\Gal(F_s/F)$, this  translates as: the field $F$ contains  a primitive $p^2$-th root of unity.] This simplifies the exposition, while  essentially making  no difference compared to the general case. However, observe that the usual restriction/corestriction argument (w.r.t. an open subgroup $H \subset G$, of prime-to-$p$ index) applies here to reduce to the  case $\F_p(1)=\F_p$-- but not to the case $\Z/p^2(1)=\Z/p^2$.
\section{Proof in the cyclotomic case, and assuming $\Z/p^2(1)=\Z/p^2$.}\label{secproof}

An essential feature of this proof, is to work integrally. For simplicity, we work over $\Z$, but for our purpose, working over $\Z_p$ would suffice.
\subsection{Resolutions in homological algebra.}\quad\\
In this section, $G$ is a finite group. The material here is classical.
\begin{defi}\label{defienG}
 Define the (augmentation) extension of $\Z[G]$-modules \[\mathcal I_G: 0 \to I_G \to \Z[G] \xrightarrow{\aug} \Z \to 0,  \]  where \[\aug([g])=1,\]  for all $g \in G$. Its kernel is the augmentation ideal of the ring $\Z[G]$.\\
For a $\Z[G]$-module $M$, one can form the extension of $\Z[G]$-modules \[\mathcal I_G \otimes_{\Z} M: 0 \to I_G \otimes M \to  M[G]\to M\to 0.  \]  Note that $M[G]=\Z[G] \otimes_\Z M$ is isomorphic to the induced module $\Ind_{\{1\}}^{G}(M)$.\\
For all $k\geq 0$,  define extensions of $\Z[G]$-modules  \[\mathcal I_G ^{k+1}:=\mathcal I_G \otimes_\Z I_G^{\otimes k}: 0 \to I_G^{\otimes (k+1)} \to I_G^{\otimes k}[G] \to I_G ^{\otimes k }  \to 0.\]  For  $n \geq 1$, concatenating the extensions $\mathcal I^ n_G, \mathcal I^{n-1}_G,\ldots, \mathcal I_G$ (in other words, forming their cup-product), one gets an $n$-extension of $\Z[G]$-modules  \[  0 \to  I_G^{\otimes n} \to  I_G^{\otimes (n-1)}[G]  \to    \ldots \to  \Z[G] \to  \Z \to 0,\] which we denote $\mathcal E^n(G)$. Denote the cohomology class of $\mathcal E^n(G)$ by \[e^n(G) \in H^n(G,I_G^{\otimes n}).\]
\end{defi}

\begin{lem}(Versality of $e^n(G)$.)\label{LemRes}\quad \\
Let $M$ be a $\Z[G]$-module. 
Then for all $n \geq 1$, the push-forward arrow
$$\fonctionnoname{\Hom_{\Z[G]-Mod}(I_G^{\otimes n},   M)}{H^n(G,M)}{f}{ f_*(e^n(G))}$$ is onto.
\end{lem}
\begin{dem}
Recall that the $\Z[G]$-modules $ I_G^{\otimes k}[G] $, for $k=0,\ldots, n-1$, are induced, or equivalently here, coinduced. Hence $$\Ext_{\Z[G]-Mod}^j(I_G^{\otimes k}[G] ,  M  )=0,$$ for all $j \geq 1$. The statement then readily follows from standard considerations in homological algebra.
\end{dem}

\subsection{The ring scheme $\RR$.}

Let $r \geq 0$ be a large enough integer. The precise meaning of `large enough'  will be clarified shortly, in Lemma \ref{lemC11vers}.   Set \[\RR:=\W_2^{[r]},\]   and denote by $(\mathbf R, \rho_\RR,\Frob_\RR, \tau_\RR)$ the associated WTF data; see \cite{F2}, section 10.
Define $q:=p^{r+1}$.
\subsection{ Hochschild cohomology of $\Z$-group schemes.}\label{sechochZ}

The goal in this section, is to build a versal $n$-extension, akin to $e^n(G)$ above, but in a more evolved algebro-geometric context- where  $\Z[G]$-modules are upgraded to  representations of smooth affine group schemes over $\Z$. For our purpose, `versal' means that the Hochschild cohomology class constructed below, should specialise to a given class in $H^n(G,\F_p)$, where $G$ is a $(1,1)$-smooth profinite group. This is made precise in  section \ref{secspe}. Eventually, in Proposition \ref{propdescent}, we show that this versal Hochschild cohomology class lifts- which implies the Smoothness Theorem by a straightforward  specialisation process. 

\begin{defi}

 Set $S:=\Spec(\Z)$. \\For an algebro-geometric structure $\mathcal X$ over $S$, denote its fiber at $p$ by \[ \overline {\mathcal X} :=\mathcal X \times_S \Spec(\F_p). \] 
\end{defi}
\begin{defi}\label{defiPk}

Let $n\geq 1$ and $r \geq 0$ be integers.

    For $k=1,\ldots,n$, introduce an extension of finite free $\Z$-modules, \[ (\mathcal P_k): 0 \to I_k \to P_k \to \Z \to 0.\]  Set \[I_{k,r+2}:=I_k \otimes_\Z\W_{r+2}(\Z),\] and \[I_{k,\RR}:=I_k \otimes_\Z\RR(\Z).\]

  \end{defi}
Next, we define natural groups that act on the data above.

  \begin{defi}
  For $k=1,\ldots,n$, consider the  $\Z$-group scheme $\Aut(\mathcal P_k),$ defined as the closed sub-group-scheme of $\GL(P_k)$ that stabilises $I_k \subset P_k$, and acts trivially on $ P_k/I_k = \Z$. In a suitable $\Z$-basis of $P_k$, it consists of matrices of the shape 
\[
\begin{pmatrix}
\ast &\ast  & \cdots & \ast & \ast \\

\vdots  & \vdots  & \ddots & \vdots  & \vdots \\
\ast  &\ast & \cdots &\ast & \ast \\
0  & 0  & \cdots & 0 & 1
\end{pmatrix}\]

Define the  green bundle  (see \cite{F2}, Definition 4.21) \[\GL(I_{k,r+2}):=\mathbf T_{r+2}(\GL(I_k)) \xrightarrow{b} \GL(I_k).\] There are natural extensions of  affine smooth $\Z$-group schemes \[1 \to \A(I_k) \to  \Aut(\mathcal P_k) \to \GL(I_k) \to 1 \]  and \[1 \to \ast \to \GL(I_{k,r+2})\xrightarrow{b} \GL(I_k) \to 1, \] 
  w.r.t. which one defines  the affine smooth $\Z$-group scheme \[ \Gamma_k:=\Aut(\mathcal P_k ) \times_{\GL(I_k)}\GL(I_{k,r+2}). \] There is a natural extension  (of smooth affine $\Z$-group schemes) \[1 \to \A(I_k) \to  \Gamma_k \to \GL(I_{k,r+2}) \to 1. \]

Define  \[ \Gamma:=\GL(I_1) \times \prod_{k=2}^n \Gamma_k,\] which fits into the  natural extension of  smooth affine $\Z$-group-schemes
 \[(\mathcal E \Gamma): 1 \to  \prod_{k=1}^n  \A(I_k) \to \Gamma \to \GL(I_1) \times \prod_{k=2}^n  \GL(I_{k,r+2}) \to 1. \]  
From now on, until specified otherwise, the acting group is $\Gamma$ (possibly upon natural group-change, e.g. via one of the projections $\Gamma \to \Gamma_k$).

\end{defi}

\begin{defi}\label{defiC}
    For $k=1,\ldots,n$, introduce the Hochschild class of the extension of $\Gamma$-bundles $(\mathcal P_k)$  (\cite{F2}, Definition 5.4), \[h(\mathcal P_k) \in  H^1(\Gamma, \A(I_k)).\] Define \[I :=\bigotimes_{k=1}^n I_k.\] The cup-product\[P:=(h(\mathcal P_1)  \cup h(\mathcal P_2)  \cup \ldots \cup h(\mathcal P_n) )\in  H^n(\Gamma, \A(I))\] is  the Hochschild class of the $n$-extension of $\Gamma$-bundles \[(\mathcal P):= 0 \to  I \to  P_1 \otimes I_2 \ldots \otimes I_n \to P_2  \otimes I_3 \ldots \otimes I_n  \to \ldots \] \[\ldots \to  P_{n-1} \otimes I_n  \to P_n \to \Z  \to 0, \] defined as the concatenation (or cup-product) of the $1$-extensions $(\mathcal P'_k)$, that are obtained by tensoring $(\mathcal P_k)$ with the $\Gamma$-bundle $ I_{k+1} \otimes \ldots \otimes I_n$. These read as\\
    
    {\centering\noindent\makebox[355pt]{$ (\mathcal {P}'_k): 0 \to  \boxed{ I_k}  \otimes I_{k+1} \otimes \ldots \otimes I_n \to \boxed{ P_k} \otimes I_{k+1} \otimes \ldots \otimes I_n \to  \boxed{\Z} \otimes I_{k+1} \otimes \ldots \otimes I_n \to 0.$}}
\end{defi}

\begin{defi}
    Define the projective bundle \[B:=\P(I) \xrightarrow{b} S.\] Observe that $B$ is a $\Gamma$-scheme in a natural way. Over $B$, denote the tautological extension of vector bundles by \[0 \to \mathcal H \to  b^*(I) \xrightarrow{\pi} \mathcal O(1) \to 0.\] Over $B$, consider the pushed-forward $n$-extension $\pi_*(b^*(\mathcal P))$.  Denote its  class by \[ C:=\pi_*(b^*(P)) \in  H^n(\Gamma, \Pi_b(\mathcal O(1))).\] Observe that $C=P$, via the natural identification $\Pi_b(\mathcal O(1)) = \A(I)$.
\end{defi}

\begin{rem}
   In this section, for  $V$  a vector bundle over $B$, the notation $V(n)$ stands for $V \otimes_{\mathcal O_B} \mathcal O(n)$. This is \textit{not} related to the   cyclotomic module. Recall that the cyclotomic module is  here the trivial module $\Z/p^2$.
\end{rem}
Think of $C$ as a universal $n$-th cohomology class with values in a line bundle- as made precise in section \ref{secspe}.

\subsection{Group-change $\Gamma_{\vers} \xrightarrow{\gamma} \Gamma$, and versal class $C_{\vers}$.}\label{secgrpch}
Recall that extensions of $\Gamma \RR$-Modules, are admissible by definition (see \cite{F2}, Definition 5.1).\\
\begin{defi}(The groups $\Gamma^1_k$, $k=2,\ldots,n$.)\label{defgamma1k}\\
Pick an integer  $k \in \{2,\ldots,n\}$.
Over $S,$ consider the extension of $\Gamma$-bundles  \[ (\mathcal P_k): 0 \to I_k \to P_k \to \Z \to 0.\]
    Denote by  \[(\mathcal E \Gamma^1_k): 1 \to \A( \Sym^q(I_k)) \to \Gamma^1_k \xrightarrow{\gamma_k} \Gamma \to 1 \] its suspension, w.r.t. the ring scheme $\RR$ (see \cite{F2}, section 13). Recall that $\gamma_k$ is the universal group-change, for the property that the extension of  $\Gamma^1_k$-bundles $\gamma_k^*(\mathcal P_k)$ lifts, to an extension of $\Gamma^1_k \RR$-Modules, denoted here by  \[ \RR(\mathcal P_{k}): 0 \to \RR(I_k) \to \ast \to \RR(\Z) \to 0.\]

\end{defi}
\begin{defi}(The group $\Gamma^1_1$.)\label{defgamma11}
Form the push-forward diagram of extension of $\Gamma$-bundles over $B$\\

{\centering\noindent\makebox[355pt]{$\xymatrix{ b^*(\mathcal P'_1): 0 \ar[r] &  I \ar[r] \ar[d]^\pi & P_1 \otimes I_2 \otimes  \ldots \otimes  I_n \ar[r] \ar[d] & I_2\otimes   \ldots \otimes  I_n \ar[r] \ar@{=}[d] & 0 \\   \mathcal P''_1:=\pi_*(b^*(\mathcal P'_1)): 0 \ar[r] &  \mathcal O(1) \ar[r] &\ast  \ar[r] &I_2 \otimes \ldots \otimes  I_n \ar[r]  & 0.}$}}
\quad\\

[In this diagram, to keep notation light, the vector bundle $b^*(I)=(I \otimes_\Z \mathcal O_B)$ over $B$, is just denoted by $I$. Same applies to other spots. This harmless abuse will be frequently committed  again, without further notice.]\\
Dualising $\mathcal P''_1$ yields the  following extension of $\Gamma$-bundles  over $B$: \[\xymatrix{(\mathcal Q_1):= 0 \ar[r] &    I_2^\vee \otimes  \ldots \otimes  I_n^\vee (1)\ar[r] &\ast  \ar[r] &  \mathcal O_B \ar[r]  & 0.}\]
    Denote by  \[(\mathcal E \Gamma^1_1):= 1 \to \A(\Sym^q( I)  \otimes \Sym^q( I_2^\vee \ldots \otimes  I_n^\vee)) \to \Gamma^1_1 \xrightarrow{\gamma_1} \Gamma \to 1 \] the suspension of  $(\mathcal Q_1)$, w.r.t. the ring scheme $\RR$. It is the universal group-change, for the property that the extension of  $\Gamma^1_1$-bundles $\gamma_1^*(\mathcal Q_1)$ lifts, to an extension of $\Gamma^1_1 \RR$-Modules  over $B$ \[\RR (\mathcal Q_1):= 0 \to \RR(I_2^\vee \otimes \ldots \otimes  I_n^\vee(1)) \to \ast \to \RR \to 0.\]

\end{defi}

\begin{defi}(The group $\Gamma_{\vers}$.)\label{defgammavers}
Define the smooth affine $\Z$-group-scheme $\Gamma_{\vers}$, as the fibered product (over $\Gamma$) of the morphisms $\gamma_1,\gamma_2, \ldots, \gamma_k$.  Denote by\\

{\centering\noindent\makebox[355pt]{$(\mathcal E \Gamma_{\vers}): 1 \to \A \left( \Sym^q( I)  \otimes \Sym^q( I_2^\vee  \otimes    \ldots \otimes  I_n^\vee)  \right )\times \prod_{k=2}^n  \A( \Sym^q(I_k)) \to \Gamma_{\vers} \xrightarrow{\gamma} \Gamma \to 1 $}}
\quad\\
the  natural extension of  smooth affine $\Z$-group-schemes.
\end{defi}

\begin{defi}

Upon group-change,  $C$ yields a class
 \[C_{\vers}:=\gamma^*( C)\in  H^n(\Gamma_{\vers},\A(I)),\] called the versal class.
\end{defi}
\subsection{Specialising $C_{\vers}$.}\label{secspe}\quad\\
First, we show that $C$ itself is versal, in the following sense.
\begin{lem}\label{lemCvers}
    Let $G$ be a profinite group. Let $D$ be a one-dimensional $(\F_p,G)$-bundle. Consider a non-zero class $c \in H^n(G,D)$. There exists $(\mathcal P_1), \ldots, (\mathcal P_n)$ as in Definition \ref{defiPk},  together with  a group homomorphism with open kernel, \[\sigma: G \to \Gamma(\Z),\] an $\F_p$-point \[ s \in B(\F_p),\] and an isomorphism of   $(\F_p,G)$-line bundles\[D \simeq \sigma^*(s^*(\mathcal O(1))),\]such that, w.r.t. this data, \[ c=\sigma^*(s^*(C)) \in  H^n(G,D).\]
\end{lem}

\begin{dem}
 Recall that \[H^n(G,D)=\varinjlim\limits_{H\trianglelefteq G }  H^n(G/H,D),\] where the  direct limit is taken w.r.t. all open normal subgroups. Therefore,  one assumes w.l.o.g. that $G$ is finite. Then, take all $(\mathcal P_k)$'s  equal to the augmentation sequence $\mathcal I_G$. The natural action of $G$ on $\mathcal I_G $ then  provides the  homomorphism $\sigma$. Applying  Lemma \ref{LemRes},  one gets a $G$-homomorphism \[f:  I=I_G ^{\otimes n}\to   D,\]  such that \[c=f_*(e^n(G)).\] Observe that $f \neq 0$ (because $c \neq 0$). Thus, by universal property of the projective space $B$, $f$ defines a point $s \in B(\F_p)$, together with a natural isomorphism  $D \simeq \sigma^*(s^*(\mathcal O(1)))$- via which the equality $c=f_*(e^n(G))$ translates as the desired statement.
\end{dem}

\begin{rem}
    Let us highlight a crucial fact, somehow hidden in the proof of the previous Lemma. There, it is essential that the group $G$ acts on $\mathcal I_G $,  and not just on the mod $p$ reduction $\mathcal I_G  \otimes_\Z \F_p$. Indeed, this ensures that the  $G$-action on  $I_G$, a priori given by a group homomorphism \[ G \to  \GL(I_G)(\Z),\] lifts through the green bundle \[\mathbf T_{r+2}(\GL(I_G)) \xrightarrow{b} \GL(I_G).\] This is because the  homomorphism induced by $b$ on $\Z$-points, reading as \[\GL(I_G)(\W_{r+2}(\Z))=\mathbf T_{r+2}(\GL(I_G))(\Z)  \xrightarrow{b} \GL(I_G)(\Z),\] has a canonical section, provided by the ring homomorphism $\Z \to \W_{r+2}(\Z)$.  This also clarifies, that $\sigma$ is well-defined. 
\end{rem}

\begin{rem}
 To define  $e^n(G)$, one considers the  `diagonal'  $n$-fold action of $G$, on  $I_G ^{\otimes n}$ (Definition \ref{defienG}). In contrast,  the class $C$  is (essentially) defined using the   action `factor by factor' of the $\Z$-group scheme $\Aut(\mathcal P_1) \times \Aut(\mathcal P_2) \times \ldots \times   \Aut(\mathcal P_n)$  on $I_1 \otimes I_2 \otimes \ldots \otimes I_n$  (Definition \ref{defienG}).  The reader may wonder, what the purpose of this `linearization' process is. As far as we can see, the answer is simple:  it renders some forthcoming computations easier, than if using the  $n$-fold action of $\Aut(\mathcal P_1)$ on  $I_1 ^{\otimes n}$, in the process of defining the versal class $C_\vers$.
\end{rem}
\begin{rem}
   Observe that the class $C$ just depends on  the ranks $r_1,\ldots,r_n$ of the free $\Z$-modules $I_1,\ldots,I_n$. In fact, for $C$ to specialise to $c$ (as Lemma \ref{lemCvers}),  these ranks may be picked all equal to $\vert G \vert-1$. This is clear from the proof above. For $m\geq 1$, take $r_i=m$ for $i=1,\ldots,n$, and  set $C_m:=C$. By what precedes, the countable family $(C_m)_{m\geq 1}$  has the property that, for any finite group $G$, for any one-dimensional $(\F_p,G)$-bundle $D$, and for any $c \in H^n(G,D)$,  one member of that family specialises to $c$ (as described in Lemma \ref{lemCvers}). This bears similarity to the  constructions made in Section 5 of \cite{DF}.
\end{rem}

\begin{lem}\label{lemC11vers}
    Let $(G,\Z/p^2)$ be a $(1,1)$-cyclotomic pair. Let $D$ be a one-dimensional $(\F_p,G)$-bundle, and let $n \geq 2$ be an integer. Consider a non-zero class $c \in H^n(G,D)$. Then, for $r$ large enough,  there exists $(\mathcal P_1), \ldots, (\mathcal P_n)$ as in Definition \ref{defiPk},  together with  a group homomorphism with open kernel, \[\sigma_{\vers}: G \to \Gamma_{\vers}(\F_p),\] an $\F_p$-point \[ s \in B(\F_p),\] and an isomorphism of   $(\F_p,G)$-bundles\[D \simeq \sigma_{\vers}^*(s^*(\mathcal O(1))),\]such that, w.r.t. this data, \[ c= \sigma_{\vers}^*(s^*(C_{\vers})) \in  H^n(G,D).\]
\end{lem}

\begin{dem}
One can first apply Lemma \ref{lemCvers}, to define $s$ and $\sigma$. It then remains to prove that  the mod $p$ reduction of $\sigma$, reading as  \[ \xymatrix{\overline \sigma: G \xrightarrow{\sigma}  \Gamma(\Z) \xrightarrow{\nat}   \Gamma(\F_p), }\] lifts to the desired $\sigma_{\vers}$.  In other words, we need to show that $\overline \sigma$ factors as \[ \xymatrix{\overline \sigma: G \xrightarrow{\sigma_{\vers}}  \Gamma_{\vers}(\F_p) \xrightarrow{\gamma}   \Gamma(\F_p). }\]

This holds for $r$ large enough: apply Proposition 14.5 of \cite{F2} (with the trivial cyclotomic twist).
\end{dem}
\quad\\

Let us proceed to proving  Theorem \ref{SmoothTh}. Let $G$ be a profinite  group, such that the pair $(G, \Z/p^2)$  is $(1,1)$-cyclotomic. Let $n \geq 2$ be an integer, and let $c \in H^n(G,\F_p)$. Let $C_{\vers}, s$ and $\sigma$ be as in 
 Lemma \ref{lemC11vers} (applied to $D:=\F_p$).\\ It then suffices to prove that $C_{\vers}$ lifts, to a class \[C_{2,\vers}\in  H^n(\Gamma_{\vers}, \Pi_b(\RR(\mathcal O(1)))).\]  Indeed, provided such a versal lift exists, simply define the specialisation \[c_2:=\sigma_{\vers}^*(s^*(C_{2,\vers}))\in  H^n(G, \Z/p^2),\] as the sought-for lift of $c$. [Recall that $\RR(\F_p)=\Z/p^2$.]\\
 
In the final steps of the proof, we prove the existence of a class $C_{2,\vers}$ as above.
\subsection{Uplifting...}\label{secUL}

\begin{defi}(Uplifting.) \\
    For $k=2,\ldots,n$, denote by \[ U_k:=\mathbf U(I_{k,r+2}) \xrightarrow{u_k} S\] the uplifting scheme of the $\W_{r+2}$-bundle $I_{k,r+2}$ (see \cite{F2}, section 12). Define the $\Gamma$-scheme \[U:=\prod_{k=2}^n  U_k \xrightarrow{u} S,\] where the product is fibered over $S$.
Define \[ T:=B \times_S U \xrightarrow{t} S.\]
\end{defi}

\begin{lem}\label{lemfil}

For $k=2,\ldots, n$, the  morphism $u_k$ is well-filtered, in a $\Gamma$-equivariant way. The good filtration of $(u_k)_*(\mathcal O_{U_k})$ is labelled by $(a_0,\ldots, a_r) \in \N^{r+1}$, well-ordered lexicographically. The corresponding graded pieces   are the $\Gamma$-vector bundles   \[\bigotimes_{i=0}^r\Sym^{a_i}   \left( \Sym^{p^i}(I_k^{(r+1-i)})\otimes \Gamma^p(\Gamma^{p^i}(I_k^{\vee(r-i)} ))\right) . \]

The   morphism $u$ is well-filtered, in a $\Gamma$-equivariant way. For our purpose, we can forget the exact shape  of the graded pieces of the associated filtration of $(u_*(\mathcal O_U)/\Z.1)$. It suffices to remember, that for each such graded piece $\mathcal Z$, there is a monomial  expression  in $(n-1)(r+2)$ vector bundles,
 \[E( V_{k,i}, \;  2 \leq k \leq n,\;   0 \leq i \leq r+1), \]  such that \[\mathcal Z=E( I_k^{(i)}, \;  2 \leq k \leq n,\;   0 \leq i \leq r+1) .\] 
 One readily checks that   its degree $(a_{k,i}) \in \Z^{(n-1)(r+2)}$ satisfies $\sum a_{k,i}=0$, and that $a_{k,i} \geq 1$ for some $(k,i) \in \{2,\ldots,n\} \times \{1,\ldots,r+1\}$.
\end{lem}

\begin{dem}
    Recall that an uplifting scheme is  a  fibered product of splitting schemes (\cite{F2}, Definition 12.16). The natural $\N$-filtrations  on  each of these  (\cite{F2}, Definition 11.2),  yield by Remark 11.5 of \cite{F2}, the sought-for good filtrations. The description of graded pieces is  slightly tedious but straightforward.
\end{dem}
\begin{lem}\label{lemextslift}

Upon base-change via  $t$ and group-change via $\gamma$, the extension \[ (\mathcal P_k): 0 \to I_k \to P_k \to \Z \to 0,\] for $k=2,\ldots,n$, acquires a natural lift,  to an extension of $\Gamma_{\vers} \RR$-bundles over $T$, \[ (\mathcal P_{k,\RR}): 0 \to I_{k,\RR} \to  P_{k,\RR} \to \RR(\mathcal O_T) \to 0.\]  Similarly, the extension  \[\xymatrix{ \mathcal P''_1: 0 \ar[r] &  \mathcal O(1) \ar[r] &\ast  \ar[r] & I_2  \otimes  \ldots \otimes  I_n \ar[r]  & 0}\] also acquires such a lift, denoted by \[ \xymatrix{(\mathcal P''_{1,\RR}): 0 \ar[r] & \RR(\mathcal O(1) )\ar[r] &\ast  \ar[r] & I_{2,\RR}  \otimes  \ldots \otimes  I_{n,\RR}\ar[r] & 0.}\]
\end{lem}

\begin{dem}
For $k=2,\ldots,n$,  apply  the Uplifting Pattern of \cite{F2}, Proposition 14.4.  Using \cite{F2}, sections 12.4 and 12.5, this Pattern  also applies to  the extension $(\mathcal Q_1)$ (see Definition \ref{defgamma11}), that thus lifts   to an extension of $\Gamma_{\vers} \RR$-bundles over $T$, \[\xymatrix{(\mathcal Q_{1,\RR}): 0 \ar[r] &    I_{2,\RR}^\vee \otimes  \ldots \otimes  I_{n,\RR}^\vee \otimes \RR(\mathcal O(1) ) \ar[r] &\ast  \ar[r] & \RR(\mathcal O_T) \ar[r]  & 0.}\] Dualising $(\mathcal Q_{1,\RR})$ yields the sought-for $(\mathcal P''_{1,\RR})$.
\end{dem}

\begin{rem}
    At the end of the previous proof, dualising is legitimate because we are dealing with extensions of $\RR$-\textit{bundles}. This is actually a major reason for introducing an uplifting scheme, in addition to   suspension, that  would be insufficient on its own. Indeed, let $V$ be a vector bundle of dimension $\geq 2$. If $r \geq 1$, it is \textit{not} true that the $\RR$-Modules $\RR(V)$ and $\RR(V^\vee)$ are dual to each other.
\end{rem}

\begin{lem}\label{lemliftoverT}
     The class  \[t^*( C_{\vers})\in  H^n(\Gamma_{\vers},\Pi_t( \mathcal O(1)))\]  has a natural lift, to a class  \[\tilde C_2 \in  H^n(\Gamma_{\vers} ,\Pi_t(\RR(\mathcal O(1)))).\] 

\end{lem}

\begin{dem}
Recall that $C$ is, by definition, the (class of the) cup-product $n$-extension \[\pi_*(\mathcal P)=\mathcal P''_1 \cup \mathcal P'_2 \cup \ldots  \cup\mathcal P'_n.\] By Lemma \ref{lemextslift}, upon the base/group-change $(t,\gamma)$, each of $\mathcal P''_1, \mathcal P'_2, \ldots, \mathcal P'_n$ lift to an extension of $\Gamma_{\vers} \RR$-bundles. The cup-product of these lifts is then an $n$-extension of $\Gamma_{\vers} \RR$-bundles , lifting the $n$-extension $\pi_*(\mathcal P)$, and reading  as\\

{\centering\noindent\makebox[355pt]{$0 \to  \RR(\mathcal O(1)) \to  \ast \to  P_{2,\RR} \otimes   I_{3,\RR}\otimes \ldots \otimes 
I_{n,\RR} \to \ldots \to    P_{n-1,\RR} \otimes I_{n,\RR}  \to P_{n,\RR} \to \RR(\mathcal O_T)  \to 0.$}}
\quad\\
Denote it by $\mathcal E_\RR$.
Consider the following  natural commutative diagram:

$$\xymatrix{ \Ext^n_{(\Gamma_\vers \RR,T)-bun  }(\RR(\mathcal O_T),\RR(\mathcal O(1)))  \ar[r]^-h \ar[d]^\rho &  H^n(\Gamma_{\vers} ,\Pi_t(\RR(\mathcal O(1))))  \ar[d]^\rho \\ \Ext^n_{(\Gamma_\vers,T)-bun  }(\mathcal O_T,\mathcal O(1))  \ar[r]^h  &  H^n(\Gamma_{\vers} ,\Pi_t(\mathcal O(1))).}$$

 This justifies, that  $\tilde C_2:=h(\mathcal E_\RR)$ does the job, as the desired lift.
\end{dem}
\quad\\
\subsection{... then descent.}\label{secdescent}
The \textit{bouquet final}   is a descent statement. The key point is to decrease the degree of the relevant cohomological obstruction, by \textit{two}: from $H^{n+1}(.)$, down  to  $H^{n-1}(.)$. This justifies most of the `integral' material developed in \cite{F2}.
\begin{prop}\label{propdescent}

   The class  \[ C_\vers \in  H^n(\Gamma_{\vers}, \Pi_b(\mathcal O(1)))\] lifts to a class  \[C_{2,\vers} \in  H^n(\Gamma_{\vers}, \Pi_b(\RR(\mathcal O(1)))).\] 
\end{prop}

\begin{dem}
The (admissible) reduction sequence of the $\RR$-bundle $ \RR(\mathcal O(1))$  reads as  \[ 0 \to \Frob_*( \mathcal O(q))  \to \RR(\mathcal O(1)) \to  \mathcal O(1) \to 0.\]

  Considering the associated exact sequence in cohomology, one sees that   the obstruction to lifting $ C_{\vers}$  is a class \[\obs (C_{2,\vers}) \in  H^{n+1}(\Gamma_{\vers},\Pi_b( \mathcal O(q))).\]
Our goal is to show \[\obs(C_{2,\vers}) \stackrel ? = 0.\] Consider the projection $\pr: T \to B$,  and the extension of $\Gamma_{\vers}$-modules over $B$, \[ 0 \to  \mathcal O_B \to \pr_*( \mathcal O_T)\to \ast \to 0. \] Upon twisting by $\mathcal O(q)$,  and using the projection formula \[\pr_*( \mathcal O_T)=u_*(\mathcal O_U)\otimes_{\Z} \mathcal O_B,\] it becomes \[ (\mathcal F): 0 \to  \mathcal O_B(q)\xrightarrow{\iota} \pr_*( \mathcal O_T(q))\to (u_*(\mathcal O_U) / \Z.1)\otimes_{\Z} \mathcal O_B(q) \to 0. \] 

By Lemma \ref{lemliftoverT}, $C_{\vers}$ lifts upon base-change to $T$. Equivalently, the corresponding obstruction vanishes over $T$. By naturality of connecting homomorphisms in cohomology, this translates as   \[\iota_*(\obs(C_{2,\vers}) )=0 \in H^{n+1}(\Gamma_{\vers},\Pi_t( \mathcal O_T(q))) . \] Considering the exact sequence in Hochschild cohomology associated to $\Pi_b(\mathcal F)$, and by dévissage along the filtration of $(u_*(\mathcal O_U) / \Z.1)$ (Lemma \ref{lemfil}), it then suffices to prove \[H^n(\Gamma_{\vers},\Pi_b(\mathcal Z(q)))  \stackrel ? = 0,\] for every graded piece $\mathcal Z$ of that filtration. By Lemma \ref{lemQFl}, such a group is of $p$-primary torsion. It thus remains to show that it is $p$-torsion-free.  Considering  (the effect in cohomology of) the extension of $\Gamma_{\vers}$-modules over $B$,  \[0 \to \mathcal Z(q) \xrightarrow{ p \Id } \mathcal Z(q) \to \overline{\mathcal Z}(q) \to 0,\] this  $p$-torsion-freeness boils down to  \[H^{n-1}(\overline {\Gamma_{\vers}},\Pi_{\overline b}(\overline{\mathcal Z}(q)))  \stackrel ? = 0,\] a fact proved in Lemma \ref{lemFp}.

\end{dem}

\begin{rem}
   The $n$-extension $(\mathcal E_\RR)$ itself does not descend to an $n$-extension of $\RR \Gamma_{\vers}$-bundles over $B$: only its Hochschild class does.
\end{rem}

\begin{lem}\label{lemQFl}

For every  graded piece $ \mathcal Z$ of the equivariant good filtration of $(u_*(\mathcal O_U)/\Z.1)$, and for every $j \geq 0,$ the group
 \[ H^j(\Gamma_{\vers}, \Pi_b( \mathcal Z(q))) \]  is of $p$-primary torsion.
\end{lem}

\begin{dem}

To begin, observe that \[ H^j(\Gamma_{\vers}, \Pi_b( \mathcal Z(q)))=  H^j(\Gamma_{\vers},  \mathcal Z \otimes_\Z \Sym^q( I)),\] by the  projection formula, combined to the classical formula $b_*(\mathcal O_B(q))=\Sym^q( I)$.
Denote by $F$  one of the following fields: $\Q$, or $\F_l$ for a prime $l \neq p$. \\ 
Proposition 4.18 of  \cite{J}, applied to $k:=\Z$ and $k':=F$,  reduces the claim to \[H^j(\Gamma_{\vers} \times_S \Spec(F), \mathcal Z \otimes_\Z \Sym^q(I) \otimes_\Z F) \stackrel {?} =0,\] for all $j \geq 0$ and all $F$ as above.
From now on,  we fix such an $F$, and work exclusively upon base-change to $\Spec(F)$. For clarity,  notation $(.)\otimes_\Z F$ and $(.)\times_S \Spec(F)$  is omitted. Observe that $\Sym^q(I)$ is a monomial expression in $I_1,I_2,\ldots, I_n$, of degree $(q,q, \ldots, q)$. By dévissage on the good filtration of Lemma \ref{lemfil},   it then suffices  to prove the following.  Let \[E \left( V_1, V_{k,i}, \;  2 \leq k \leq n,\;   0 \leq i \leq r+1) \right)\] be a monomial expression in $1+(n-1)(r+2)$ variables, whose degree w.r.t. $ V_{k,i}$ is non-zero, for some  $k \in \{2,\ldots,n\}$ and $i \in \{1,\ldots,r+1\}$. Then one has \[ H^j(\Gamma_{\vers}, E(I_1, I_k^{(i)}))\stackrel{?} =0. \] 
Recall that, since $p \in F^\times$,  the ring scheme $\W_{r+2}$ is isomorphic, over $F$, to the direct product $\W_1^{r+2}$. As opposed to what happens in characteristic $p$,  here `Frobenius twists are separated':  for $k=2,\ldots,n$,  \[I_{k,r+2}=\underbrace {I_k \bigoplus I_k^{(1)}  \bigoplus \ldots  \bigoplus  I_k^{(r+1)}}_{r+2 \; \mathrm{factors}}\] and compatibly \[ \GL(I_{k,r+2})=\underbrace {\GL(I_k) \times \GL(I_k^{(1)}) \times\ldots \times \GL(I_k^{(r+1)})}_{r+2 \; \mathrm{factors}}.\]

 Define\[\mathbf G:=  \GL(I_1) \times \prod_{k=2}^n \GL(I_{k,r+2})= \GL(I_1) \times \prod_{\substack{2\leq k\leq n\\ 0\leq i \leq r+1}}\GL(I_k^{(i)}).\] Note that, up to isomorphism,  
 $\mathbf G$ is just a finite product of $\GL_d$'s, for various integers $d$. Consider the group extensions $(\mathcal E \Gamma_{\vers})$ and $(\mathcal E \Gamma)$. Their kernels share a common shape: indeed, both are products of (affine spaces of) monomial expressions  in $I_1,\ldots,I_n$. Observe that \textit {these  do not involve any $I_k^{(i)}$ for $i \geq 1$,}  and that they act trivially on the coefficients of the cohomology. \noindent 

Using the spectral sequence of the group extension $(\mathcal E \Gamma_{\vers})$, then that of  $(\mathcal E \Gamma)$, and  each time applying Lemma \ref{lemSS},  one reduces to the following vanishing statement. For every monomial expression \[E \left( V_1, V_{k,i}, \;  2 \leq k \leq n,\;   0 \leq i \leq r+1) \right)\] with the properties  previously stated,  and for every $j \geq 0$, one has
 \[ H^j(\mathbf G,   E(I_1, I_k^{(i)}))\stackrel{?} =0.\]
 
This is proved in \cite{F2}, Lemma 7.3.\end{dem}
\quad\\

\begin{lem}\label{lemFp}
For every graded piece $\mathcal Z$ of the good filtration of $u_*(\mathcal O_U)$, 
    the group \[ H^{n-1}(\overline {\Gamma_{\vers}}, \Pi_{\overline b}(\overline { \mathcal Z} \otimes \mathcal O(q)))\] vanishes.
\end{lem}

\begin{dem}
In this proof, we work exclusively  upon base-change to $\Spec(\F_p)$. For clarity, \textit{overlining is omitted from notation}. The strategy is similar to that of Lemma \ref{lemQFl}, of which we keep notation.  Here also, the desired vanishing statement is equivalent to  that of \[H^{n-1}(\Gamma_{\vers}, \mathcal Z \otimes \Sym^q(I)) . \] There is one significant difference: working here mod $p$, `frobenius twists are inseparable': in characteristic $p$, the frobenius $(V \mapsto V^{(1)})$ is a monomial functor,  of degree $p$. Consequently, the monomial expressions giving the graded pieces of the good filtration of $u_*(\mathcal O_U)$ in Lemma   \ref{lemfil}, are `simpler':  they are of the shape $E(I_2, \ldots, I_n)$,  where   $E(V_2, \ldots, V_n)$ is a  monomial expression, of degree zero.  Taking into account that the degree of  $\Sym^q(I)$ is  $(q,q, \ldots,q)$, the statement boils down, by dévissage, to \[H^{n-1}(\Gamma_{\vers}, E(I_1, I_2, \ldots, I_n))  \stackrel {?} = 0, \] for every  monomial expression $E(V_1, V_2, \ldots, V_n)$, of degree $(q,q, \ldots,q)$. Also, working here over $\F_p$, the group \[\mathbf G:=  \GL(I_1) \times \prod_{k=2}^n \GL(I_{k,r+2})\]   has a natural composition  series, whose top graded piece  is \[\prod_{\substack{1\leq k\leq n}} \GL(I_k),\] and other graded pieces are the vector groups \[\A(\End(I_k)^{(i)}), \]for $2\leq k\leq n$ and $ 0\leq i\leq r+1$.
    Clearly, the latter are   monomial expressions $E(I_2,\ldots,I_n)$,  homogeneous of degree  zero.  
   The kernels of $(\mathcal E \Gamma)$  and  $(\mathcal E \Gamma_{\vers})$ are the same monomial expressions, as in the proof of Lemma \ref{lemQFl}: they do not involve frobenius twists, and  they  act trivially on the coefficients of the cohomology.  We shall need here an extra sip of information: the degree of such an expression is $(0,\ldots,0,a,0,\ldots,0)$, where $a=1$ (case of $(\mathcal E \Gamma)$) or $a=q$ (case of $(\mathcal E \Gamma_{\vers})$), and $a$ occurs as the $k$-th entry, for some $k \in \{1,\ldots,n\}$.\\
Altogether, one gets a composition series (of linear algebraic $\F_p$-groups) \[1= \mathbf  K_0 \subset \mathbf  K_1  \subset  \ldots \subset \mathbf  K_N = \Gamma_{\vers},\] such that \[ \mathbf  K_N/  \mathbf  K_{N-1}=\prod_{\substack{1\leq k\leq n}} \GL(I_k),\] and all other graded pieces $\mathbf  K_{j+1}/ \mathbf K_j$ are   monomial expressions of the shape $E_j(I_1,I_2,\ldots,I_n)$, whose degree   is  $(0,\ldots,0,+,0,\ldots,0)$, where $+$ denotes a  non-negative integer.
 To conclude,  let us show the following key statement. \\
  Let $m \in \{0 ,\ldots,  n-1\}$, and let $a_1,\ldots, a_n$  be non-negative integers, such that at most $m$ of them are non-zero. Let  $ E(V_1,\ldots,V_n)$ be a monomial expression,  of degree  $ (q-a_1,q-a_2,\ldots,q-a_n).$ Then, for every $j \in \{0, \ldots, N-1\}$, \[ H^{n-1-m}\left(\Gamma_\vers/ \mathbf K_j, E(I_1,I_2,\ldots, I_n)\right) \stackrel ? = 0.\]

The proof of this statement is by descending induction on  $j$. The starting case $j=N-1$ holds by \cite{F2}, Lemma 7.3: indeed, the assumptions imply that the $n$ integers $(q-a_i)$ cannot all vanish. The induction step (going from $j+1$ to $j$) is  done using the  spectral sequence of  the group extension  \[1 \to E_j(I_1,I_2,\ldots,I_n) \to\Gamma_\vers/ \mathbf K_j  \to  \Gamma_\vers/ \mathbf K_{j+1} \to 1.\]
 Indeed, applying Lemma \ref{lemSS} to this extension (over $k=\F_p$), it suffices to prove 
\[ H^{n-1-m}\left(\Gamma_\vers/ \mathbf K_{j+1}, E(I_1,I_2,\ldots, I_n)\right) \stackrel ? = 0,\] which holds by the induction hypothesis, and \[ H^{n-1-m'}\left(\Gamma_\vers/ \mathbf K_{j+1}, E'(I_1,I_2,\ldots, I_n)\right) \stackrel ? = 0,\]  for all $m <m' \leq n-1$, and for every  monomial expression $E'(.)$, whose degree reads as  $ (q-a'_1,q-a'_2,\ldots,q-a'_n),$ where the  integers $a'_i \geq 0$, are such that at most $(m+1)$ of them are non-zero. [Precisely: $a'_i=a_i$ for all indices $i$, but possibly one.] Since $m+1 \leq m'$, this case is also covered by the induction hypothesis.\\
The statement thus holds for $j=0$ and $m=0$, yielding the desired vanishing.

\end{dem}

\begin{lem}\label{lemSS}
 Let $k$ be a perfect field. Let $V$ be a finite-dimensional $k$-vector space. Consider an extension of linear algebraic $k$-groups \[(\mathcal E): 1 \to \A(V) \to \mathbf H^1 \to \mathbf H \to 1.\] Let  $C$  be a $k$-representation of $\mathbf H^1$, restricting to a trivial representation of $\A(V)$. Consider the spectral sequence relative to the group extension $(\mathcal E)$, \[H^i(\mathbf H, H^j( \A(V) ,C)) \Rightarrow H^{i+j}(  \mathbf H^1 ,C).\] Then, the left side is a direct sum of (possibly infinitely many) groups of the shape
 \[H^i(\mathbf H,E(V,C)),\] for  monomial expressions $E(V,C)$  of degree $(-d,1)$, for (various) $d \geq 1$.
\end{lem}

\begin{dem}
   The groups  $H^j( \A(V) ,C)$ are computed in  \cite{J}, section 4.27, for $C=k$. For the general case, just apply flat base-change:\[H^j( \A(V) ,C)= H^j( \A(V) ,k) \otimes_k C.\] The result depends on the characteristic of $k$ (zero, $p=2$ or $p$ odd). However, in all three cases, it is straighforward to see, that it is a  direct sum of monomial expressions $E(V,C)$ of the required form. To conclude, use that $H^i(\mathbf H,.)$ commutes to direct limits (\cite{J}, Lemma 4.17).
\end{dem} 
\section{Proof of the Smoothness Theorem, general  case.}\label{secproof1}
Section \ref{secproof} adapts, to deal with an arbitrary  $(1,1)$-cyclotomic pair $(G,\Z/p^2(1))$. \\We will  browse through the proof, highlighting the steps that  need to be modified.\\ The reader is assumed to be familiar with  section 14.2 of \cite{F2}, especially Proposition 14.10, Lemma 14.11 and Proposition 14.12 thereof. In particular, recall that  \[\Z/p^2[1]:=\Z/p^2(1) \otimes_\Z \W_2(\F_p(-1))\] is another cyclotomic module for $G$, such that $\F_p[1]=\F_p$.\\
Until the end of section \ref{sechochZ}, there is nothing to change. Modifications begin in Definitions \ref{defgamma1k}, \ref{defgamma11} and \ref{defgammavers}, where one incorporates a scheme-theoretic cyclotomic twist, as follows. In place of $\Gamma$, one uses \[\Gamma_1:=\Gamma \times_{\Z} \G_{a/m}, \] and the related $\Gamma_1 \RR$-bundle $\RR(\Z)[1]$, over $\Spec(\Z)$. Using the suspension process of Definition 14.9 of \cite{F2}, the modified group-changes then look like

  \[(\mathcal E \Gamma^1_k): 1 \to \A( \Sym^q(I_k) \{q\}) \to \Gamma^1_k \xrightarrow{\gamma_k} \Gamma \times_{\Z} \G_{a/m}\to 1, \] 

   \[(\mathcal E \Gamma^1_1):= 1 \to \A(\Sym^q( I)  \otimes \Sym^q( I_2^\vee \ldots \otimes  I_n^\vee)\{q\}) \to \Gamma^1_1 \xrightarrow{\gamma_1} \Gamma \times_{\Z} \G_{a/m}\to 1, \]
\quad\\   
and \\

{\centering\noindent\makebox[355pt]{$(\mathcal E \Gamma_{\vers}): 1 \to \A \left( \Sym^q( I)  \otimes \Sym^q( I_2^\vee  \otimes    \ldots \otimes  I_n^\vee) \{q\} \right )\times \prod_{k=2}^n  \A( \Sym^q(I_k)\{q\}) \to \Gamma_{\vers} \xrightarrow{\gamma} \Gamma \times_{\Z} \G_{a/m} \to 1 .$}\\\quad\\}

For $k=2, \ldots,n$, upon group-change to $\Gamma_\vers$, $\mathcal P_k$ lifts  to an extension of $\Gamma_{\vers} \RR$-Modules  \[ 0 \to \RR(I_k)[1] \to \ast \to \RR(\Z) \to 0.\] Similarly,  $\mathcal Q_1$ lifts, to an extension of $\Gamma_{\vers} \RR$-Modules  over $B$ \[0 \to \RR(I_2^\vee \otimes \ldots \otimes  I_n^\vee(1))[1] \to \ast \to \RR \to 0.\]

[Note the presence of two twists here. These are of very different origins.]\\
Here again, the versal class is
 \[C_{\vers}:=\gamma^*( C)\in  H^n(\Gamma_{\vers},\A(I)).\]

 The specialisation process (Lemmas \ref{lemCvers} and \ref{lemC11vers}) is the same,  using Proposition 14.12 of \cite{F2}, instead of 14.5. In section \ref{secUL}, the uplifting schemes $u_k$ and $u$ remain unchanged  (these have nothing to do with  the cyclotomic twist $.[1]$).  The versal lift in Lemma \ref{lemliftoverT} (whose proof is the same) reads here as  \[\tilde C_2 \in  H^n(\Gamma_{\vers} ,\Pi_t(\RR(\mathcal O(1)))[n]).\]  Observe that the correct cyclotomic twist is indeed \[.[n]:=. \underbrace{\otimes [1] \otimes \ldots \otimes [1]}_{n \; \mathrm{times}},\]  because $\tilde C_2$ is the cup-product of the  extensions (of $\Gamma_\vers \RR$-bundles over $T$)
 \[ (\mathcal P_{k,\RR}): 0 \to I_{k,\RR}[1] \to  P_{k,\RR} \to \RR(\mathcal O_T) \to 0,\]resp.
 \[ \xymatrix{(\mathcal P''_{1,\RR}): 0 \ar[r] & \RR(\mathcal O(1) )[1]\ar[r] &\ast  \ar[r] & I_{2,\RR}  \otimes  \ldots \otimes  I_{n,\RR}\ar[r] & 0,}\] that lift $  (\mathcal P_k)$ ($2 \leq k \leq n$), resp. $(\mathcal P''_1)$.\\
 It remains to explain how to adapt the descent statement of section \ref{secdescent}: Proposition \ref{propdescent}. The relevant admissible extension at the beginning of its proof, is

  \[ 0 \to \Frob_*( \mathcal O(q)\{q\})  \to \RR(\mathcal O(1))[1] \to  \mathcal O(1) \to 0.\]
The same proof goes through, provided one can show that the twisted versions of Lemma \ref{lemQFl} hold; namely, that the groups

\[ H^j(\Gamma_{\vers}, \Pi_b( \mathcal Z(q))\{q\}) \] are of $p$-primary torsion, and that the groups  \[ H^{n-1}(\overline {\Gamma_{\vers}}, \Pi_{\overline b}(\overline { \mathcal Z} \otimes \mathcal O(q)))\] vanish (observe that the twist $\{q\}$ disappears upon reduction mod $p$).\\
Recall the (slightly modified) definition of the versal `twisted' group $\Gamma_{\vers}$ considered here. For the vanishing of the former groups, the argument is almost the same- one just adds an extra factor $(. \times \G_m)$ to the group $\mathbf G$, corresponding to (the fiber over $F$  of) $\G_{a/m}$. For the latter groups, recall that the fiber of $\G_{a/m}$  over  $\F_p$ is $\G_a$. Accordingly, the  proof of Lemma \ref{lemFp} adapts to the twisted setting- taking into account  one extra term in the  composition series\[1= \mathbf  K_0 \subset \mathbf  K_1  \subset  \ldots \subset \mathbf  K_N = \Gamma_{\vers},\] that contributes to a graded piece $\G_a$,  acting trivially. One readily checks, that this does not affect the proof of the `key statement' thereafter. This concludes the proof of the smoothness theorem, in the $(1,1)$-cyclotomic case.\\
  
The case of a $(1,1)$-smooth profinite group $G$ is then straightforward. Indeed, by Theorem 12.6 of \cite{DCF1}, the group $G$ is $(1,1)$-cyclothymic. One may then use the `cyclothymic' variant of Theorem A of  \cite{DCF1}: Theorem 12.5 thereof. Indeed, observe that Theorem A (a.k.a. Proposition 8.2 of \cite{F2}), in the  $(1,1)$-cyclotomic case, is actually applied just once: in the specialisation process, when invoking Proposition 14.12 of \cite{F2}. Therefore, in  the $(1,1)$-smooth case, one can use a `local' cyclotomic module $\Z/p^2[c]$, in place of a `global' $\Z/p^2(1)$. Accordingly, in Proposition 14.12 of \cite{F2}, the character $\chi$ (used for specialisation) is that of  $\Z/p^2[c]$.

\section{What next?}

We suggest  topics of investigation, as a continuation of this work.
\subsection{The Symbols conjecture}
Roughly speaking, the following Question tentatively describes an essential trait  of $(1,1)$-smooth profinite groups.

\begin{qu}
     Let $\mathcal P(G)$ be a property, depending on a profinite group $G$. Assume that $\mathcal P(G)$ holds true for all groups of the shape $G=\Gal(E/F)$, where the (possibly infinite) Galois field extension $E/F$ is such that these conditions simultaneously hold. \begin{enumerate}
     
     \item {The extension $E/F$ is a pro-$p$-extension, of fields of characteristic $\neq p$.} \item{The field $E$ contains all $p$-th roots of unity.}\item{The $p$-power map $E^\times \xrightarrow{a \mapsto a^p} E^\times$  is onto.}
 \end{enumerate}
Does $\mathcal P(G)$  hold true for all  $(1,1)$-smooth pro-$p$-groups $G$?
\end{qu}

 \noindent We do not expect the answer to be  positive for all properties $\mathcal P$. It would be interesting to characterise a nice class of these, for which it  actually is.
 
 \begin{rem}
 Recall that conditions (2) and (3) hold, if $E=F_{sep}$ (a separable closure of $F$). They are sufficient to ensure that $\Gal(E/F)$ is  $(1,1)$-smooth.
 \end{rem}
 
We now make precise a conjecture (or question, depending on the reader's taste).

\begin{defi}[{$\infty$-smooth groups}]\quad\\
A profinite group $G$ is $\infty$-smooth if the following holds. Consider a finite collection $H_1,\ldots,H_N \subset G$ of open subgroups of $G$.\\
For each $i=1,\ldots,N$, let $$c_i \in H^{n_i}(H_i,\F_p)$$ be a cohomology class, for various positive integers $n_i$'s. Put $C=(c_1,\ldots,c_N)$.\\
Then, there exists a continuous character $$\chi(C): G \to \Z_p^\times,$$ with trivial mod $p$ reduction $G \to \F_p^\times$, such that the following holds.\\
Denote by $\Z_p(C)$ the group $\Z_p$, on which $G$ continuously acts via $\chi(C)$.\\
Then, for each $i=1,\ldots,N$, the class $c_i$ lifts through the natural map $$H^{n_i}(H_i,\Z_p(C)^{\otimes n_i}) \to H^n(H_i,\F_p).$$ The $G$-module $\Z_p(C)$  is called `cyclotomic w.r.t. $c_1,\ldots,c_N$'.
\end{defi}
\begin{rem}
If $(G,\Z_p(1))$ is an $(n, \infty)$-cyclotomic pair for all $n \geq 1$, then the group $G$ is $\infty$-smooth: one can choose  $\Z_p(C):=\Z_p(1)$ in  Definition above. Roughly speaking, $\infty$-smooth profinite groups correspond to such pairs, but without  a `global' cyclotomic module.
\end{rem}

\begin{defi}
    Let $H$ be a profinite group. Let $x_1, \ldots, x_n \in H^1(H,\F_p)$. The element \[ x_1 \cup \ldots \cup x_n \in H^n(H,\F_p)\] is called a symbol (of degree $n$). If $H$ is an open subgroup of another profinite group $G$, the element  \[ \Cor_H^G(x_1 \cup \ldots \cup x_n) \in H^n(G,\F_p)\] is called a quasi-symbol (of degree $n$).
\end{defi}
\begin{defi}
    Let $V$ be an $\F_p$-vector space. Denote by \[T^*V:=\bigoplus_{n\geq 0}T^n(V)=\bigoplus_{n\geq 0} \underbrace{V \otimes \ldots \otimes V}_{n \; \mathrm{times}}\] its tensor algebra.
\end{defi}
 Let us mention that items (2) and (3) of the following Conjecture,  are  currently being investigated by other researchers, independently and under a different perspective. We refer to their preprint \cite{CMOSUB}, where among other noteworthy results, the authors determine all finite groups $G$ that satisfy item (2).

\begin{conj}[Symbols Conjecture]\label{SymbolsConj}\hfill \\
Let $G$ be a $(1,1)$-smooth pro-$p$-group.  Then the following holds.

\begin{enumerate}
\item{The group $G$ is  $\infty$-smooth.}
       \item{(Degree one generates the mod $p$ cohomology algebra) \\For all $n \geq 2$, the cup-product \[T^n H^1(G,\F_p) \to  H^n(G,\F_p),\] \[x_1 \otimes \ldots \otimes x_n \mapsto x_1 \cup \ldots \cup x_n \] is surjective. Equivalently, $H^n(G,\F_p)$ is additively generated  by    symbols.} \item{(Strengthening of (2): relations are purely quadratic).\\ The kernel of the surjection (of graded rings) \[T^* H^1(G,\F_p) \to H^*(G,\F_p)\] is the ideal generated by  pure tensors of degre $2$, that are of the shape \[ x \cup \Cor_x(y)\in H^2(G,\F_p),\] where $x \in H^1(G,\F_p)=\Hom(G,\F_p)$, where $y \in H^1(\Ker(x),\F_p)$, and where \[ \Cor_x:H^1(\Ker(x),\F_p) \to H^1(G,\F_p) \] denotes the corestriction (=norm), w.r.t. the subgroup $\Ker(x) \subset G$.\\ \noindent [If $x\neq 0$, it is a normal subgroup of index $p$.]}
\end{enumerate}
\end{conj}

\begin{rem}
The description of relations appearing in (3), is reminiscent of  the Steinberg relation $(x) \cup (1-x)=0$, in Milnor $K$-theory of a field $F$ ($x \in F-\{0,1\}$). Arguably, (3) is what comes closest to this, in the realm of group cohomology.
\end{rem}

\begin{rem}
    Let $G$ be an absolute Galois group.  Let $n\geq 1$ be an integer. Let $x \in H^1(G,\F_p)$ and $y \in H^n(G,\F_p)$. Using the norm-residue isomorphism theorem, it is classical that the following  are equivalent: \begin{enumerate}
        \item{$x \cup y=0 \in H^{n+1}(G,\F_p)$,} \item{$y$ is in the image of $ \Cor_x:H^n(\Ker(x),\F_p) \to H^n(G,\F_p) . $}
    \end{enumerate} Note that the implication $ (2) \Rightarrow (1)$ holds for every profinite group $G$, and that so does the converse if $p=2$. With the help of theorem \ref{SmoothTh},  it is an interesting (and challenging) exercise to prove that the equivalence actually holds for every $(1,1)$-smooth profinite group.
\end{rem}
\begin{rem}
Point (2) of the Symbols Conjecture is (obviously) equivalent to the combination of the following two statements.
\begin{enumerate}[(2a)]
    \item  Quasi-symbols are sums of symbols.
    \item For all $n \geq 2$,  $H^n(G,\F_p)$  is generated by quasi-symbols.
\end{enumerate}

For $(1,1)$-smooth profinite groups, our intuition is that (2b)  holds. However, we express a serious reservation about (2a), for the following reason. In Galois cohomology, the usual proof that (2a) holds does not at all use Kummer theory. Instead, it relies on the existence of the Galois symbol, from Milnor $K$-theory to Galois cohomology, and on the Bass-Tate Lemma in Milnor $K$-theory. In turn, the effective proof of this Lemma crucially exploits the  Euclidean division algorithm for polynomials  (for details, see \cite{GS}, especially 7.2.) Thus,  it may well be the case, that (2a) is of a whole other nature than Kummer theory- and even more so, if one relaxes the assumption that $G$ is pro-$p$. This uncharted mathematical territory, is still mysterious for us.
\end{rem}

\begin{exo}
Show in an elementary way, that property (2b) above, implies the Norm Residue Isomorphism Theorem.
\end{exo}

\section*{Appendix: more cyclotomic pairs.}
\subsection{Demushkin groups.}\quad\\
Let $G$ be a Demushkin group; that is: a pro-$p$-group such that $H^2(G,\F_p) \simeq \F_p$, and satisfying Poincaré duality in dimension 2. \\Then, by \cite{Se2}, 9.3, there exists a cyclotomic module $\Z_p(1)$, such that the pair $(G,\Z_p(1))$ is $(1, \infty)$-cyclotomic.  [For an explicit and slightly stronger result,  see \cite{CDF}, Proposition 5.1.] Observe that the validity of the Smoothness Theorem, in this case, is an easy fact: indeed, by \cite{Se2}, Demushkin groups are of cohomological dimension 2- except in the (trivial) case $p=2$ and $G=\Z/2$.

\subsection{Curves over algebraically closed fields.}\quad\\
In this section, $F$ is an algebraically closed field.\\
Let $C$ be an $F$-curve.  [= a one-dimensional separated $F$-scheme of finite-type.] To fix ideas, assume that $C$ is connected. Denote `its' étale fundamental group by $G:=\pi_1(C)$. If $C/F$ is smooth, then  $(G,\Z/p^2)$ is $(1,1)$-cyclotomic, by \cite{DCF0}, Proposition 4.11.  As shown next, the smoothness assumption is in fact superfluous.

\begin{defi}
Let $C$ be a connected and reduced $F$-curve.  Denote by $$\Pic^0(C) \subset \Pic(C)$$ the subgroup consisting of isomorphism classes of line bundles on $C$, whose degrees w.r.t. all proper irreducible components of $C$ are $0 \in \Z$. \\
{[}Note that  $\Pic^0(C) = \Pic(C)$, if and only if $C$ is affine.]
\end{defi}
\begin{prop}
Let $C$ be a connected $F$-curve.  Denote its étale fundamental group by \[G:=\pi_1(C).\]   The following holds.

\begin{enumerate}
    \item{Assume  that $C$ is reduced. Let $\overline C$ be a proper  $F$-curve, containing $C$ as an open subvariety, and such that  $\overline C \backslash C$ consists of (finitely many) smooth $F$-points. Then, the restriction map $\Pic^0(\overline C) \to \Pic^0(C)$ is surjective.} \item{If $F$ has characteristic not $p$, the pair $(G,\Z_p)$ is $(1,\infty)$-cyclotomic.  }  \item{If $F$ has characteristic $p$, the pair $(G,\Z_p)$ is $(1,\infty)$-cyclotomic. }
\end{enumerate}
\end{prop}

\begin{dem}
Let $x \in  \Pic^0(C)$. Let us show that $x$ extends to an element of $\Pic^0(\overline C)$. Denote by $C_{sm} \subset C$ the smooth locus of $C$. Since $C$ is reduced, $C_{sm}$ is the complement of a finite set of closed points. We proceed by induction, on the smallest cardinality of a finite subset $Q \subset C_{sm}(F)$, such that $\Res(x)=0 \in  \Pic^0(C-Q)$. Clearly, the result holds when $Q$ is empty.

For the induction step, pick $q \in Q$. The class $\Res(x) \in  \Pic^0(C-\{q\})$ extends to  $\overline x \in  \Pic^0(\overline C)$, by induction (applied to $C-\{q\}$). Changing $x$ to $(x- \overline x)$, we are reduced to the case $\Res(x)=0 \in  \Pic^0(C-\{q\})$. In other words, $x=\mathcal O(aq)$ for some $a \in \Z$. Let  $C_q$ be the irreductible component of $C$  which contains $q$. If $C_q$ is proper, taking degree  w.r.t. $C_q$ yields $a=0$, hence $x=0$ and we are done. If $C_q$ is affine, there exists, by assumption, a smooth point $\overline c_q \in \overline C_q(F)-C_q(F)$. Then  $$(\mathcal O(aq)- \mathcal O(a\overline c_q)) \in \Pic^0(\overline C)$$ is the sought-for extension of $x$.

\noindent To prove item (2), we can replace $C$ by its maximal reduced subscheme $C_{red}$. Henceforth, we assume that $C$ is reduced. Assume first, that $C$ is proper. We then know that $\Pic^0_{C/F}$ is an extension of an Abelian variety $A$ by a commutative smooth connected linear algebraic group $L$ (see \cite[Chap. 9, Cor. 11]{BLR}). Since $p \neq 0 \in F$,  both $A(F)$ and $L(F)$ are $p$-divisible (in the naive sense, that their $p$-th power map $(u \mapsto pu)$ is onto). Hence, so is the group $ \Pic^0_{C/F}(F)$. Moreover, the Néron-Severi group of $C$ is torsion-free (loc. cit., Corollary 14). Using these two facts, the proof of \cite{DCF0}, Proposition 4.11, adapts verbatim. To deal with the general case, where $C/F$ is not assumed to be proper,  consider an $F$-compactification $C \subset \overline C$ as in item (1). The same proof as that given in  \textit{loc. cit.} then goes through.\\
Item (3): the proof for smooth curves  in \textit{loc. cit.} actually works in general.\end{dem}

\bibliographystyle{plain}
\bibliography{biblitex.bib}

\end{document}